\documentclass[reqno]{amsart}

\usepackage[centertags]{amsmath}

\usepackage{newlfont}

\usepackage{latexsym,enumerate,epsfig,listings}
\usepackage{amsthm,amsopn,amstext,amscd,amsfonts,amssymb}

\usepackage{graphicx,psfrag}

\newcommand{\beeq}{\begin{equation}}
\newcommand{\eneq}{\end{equation}}
\newcommand{\bearno}{\begin{eqnarray*}}
\newcommand{\enarno}{\end{eqnarray*}}
\newcommand{\befi}{\begin{figure}}
\newcommand{\enfi}{\end{figure}}

\def\ep{\hfill $\Box$}

\def\bp{\noindent{\bf Proof.}\ }

\def\R{\mathbb R}
\def\I{{\cal I}}
\def\C{{\cal C}}
\def\D{{\cal D}}

\def\Z{\mathbb Z}

\def\A{{\cal A}}
\def\B{{\cal B}}

\def\SS{{\cal S}}
\def\S{{\cal S}}

\newtheorem{remark}{Remark}[section]

\newtheorem{example}{Example}
\newtheorem{theorem}{Theorem}
\newtheorem{proposition}{Proposition}

\begin{document}

\title[Stability of birth-and-death processes]{Stability of multi-dimensional birth-and-death processes with state-dependent $0$-homogeneous jumps}

\author[M. Jonckheere and S. Shneer]{M. Jonckheere and S. Shneer \\ \vspace{.07in} \\
{\it Department of Mathematics and Computer Science \\ Eindhoven University of Technology }\\{\tiny P.O. Box 513, 5600 MB \\ Eindhoven, The Netherlands }}

\address{Matthieu Jonckheere \\ Department of Mathematics and Computer Science \\ Eindhoven University of Technology \\ P.O. Box 513, 5600 MB \\ Eindhoven, The Netherlands \\ E-mail: {\tt m.t.s.jonckheere@tue.nl}}
\address{Seva Shneer \\ Department of Mathematics and Computer Science \\ Eindhoven University of Technology \\ P.O. Box 513, 5600 MB \\ Eindhoven, The Netherlands \\ E-mail: {\tt shneer@eurandom.tue.nl}}

\maketitle

\begin{abstract}

We study the positive recurrence of multi-dimensional birth-and-death processes
describing the evolution of a large class of stochastic systems, a typical example being
the randomly varying number of flow-level transfers in a telecommunication wire-line or wireless network. 

We first provide a generic method to construct a Lyapunov function when the drift can be extended to a smooth function on $\mathbb R^N$,
using an associated deterministic dynamical system. This approach gives an elementary proof of ergodicity without needing to establish the convergence of the scaled version of the process towards a fluid limit and then proving that the stability of the fluid limit implies the stability of the process. We also provide a counterpart result proving instability conditions.

We then show how discontinuous drifts change the nature of the stability conditions
and we provide generic sufficient stability conditions
having a simple geometric interpretation. These conditions turn out to be necessary (outside a negligible set of the parameter space) for piece-wise constant drifts in dimension $2$.

\end{abstract}

\section{Introduction}

We study the stochastic stability of multi-dimensional birth and death processes $X=(X_1,\ldots,X_N)$ on $\mathbb Z_+ ^N$, ($N$ being an integer greater than 1) with state-dependent birth and death rates (respectively $\lambda(x)=(\lambda_i(x))_{i=1\ldots N}$ and $(\phi_i(x))_{i=1\ldots N}$, with $x=(x_1,\ldots,x_N)$) being $0$-homogeneous functions, i.e. such that
$\lambda( \alpha x)=\lambda(x)$ for any $\alpha >0$ and for any $x \in \mathbb Z_+^N$.

This assumption arises naturally in many queueing networks representing communication
or manufacturing systems where the server capacity depends on the states of all queues. Cellular radio networks are a typical example: available transmission rate for customers in a particular cell is decreasing when the number of customers in the neighboring cells increases~\cite{bonald2004a}.
More generally, many telecommunication networks can be represented (at a sufficiently large time scale) as processor sharing
networks, 
with processing rates that may depend on the number of customers
at each node of the network~\cite{bonald2006}.
The notion of fairness between the classes of customers has led to the introduction of allocations of bandwidth naturally depending on the proportion of customers of each class rather than on the number of customers \cite{massoulie2007}. This justifies our assumption of $0$-homogeneity.

The assumptions made here are also relevant for describing load balancing schemes between a set of servers
or computer systems. Very few results have been discussed in the literature on simple schemes such
as joining the shortest queue when the death rates are not constant. (For the constant case, see \cite{foley2001}).

\

A general framework for analyzing stochastic stability
consists in applying the Foster--Lyapunov criteria, which are based on finding
a suitable test function having a positive or negative mean drift
in almost all states of the state space~\cite{tweedie1978,fayolle1995,meyn1993}.
When there are further restrictive assumptions made on service rates, an appropriate Lyapunov function has been found in many cases. This is, for example, possible for rates being the solution of specific optimization problems \cite{bonald2006},
or for small dimensions (2 and 3) when the rates of the process are constant on sub-faces of the orthant \cite{fayolle1995}.
For more complex systems, however, finding a Lyapunov function can be a formidable task.

An alternative tool for deriving stability conditions is to study whether the system of interest
is stochastically comparable to another system that is easier to analyze. This approach
was first used in the multi-class queueing context by Rao and Ephremides~\cite{rao1988} and Szpankowski~\cite{szpankowski1988},
and later refined by Szpankowski~\cite{szpankowski1994}, to characterize the stability of buffered random access systems.
It was later generalized to birth-and-death processes with state-dependent transitions with fixed birth rates and decreasing death rates with uniform limits in \cite{borst2008}. This approach relies however on quite specific assumptions that are not verified even for simple processor sharing systems.

Finally, many stability results have been obtained using the so-called ODE (ordinary differential equations)
methods. A powerful exposition of these ideas applied to controlled random walks can be found in Chapter 10 of \cite{meyn2007} and in \cite{fort2008}. 
The use of ODE methods is usually coupled with the analysis of fluid limits:
first the convergence of a scaled version of the process towards a fluid limit is proven; then
(under restrictive conditions)
the stability of the fluid limit is proven to imply the stability (positive recurrence in our case)
of the stochastic process.
Stability conditions for a wide class of multi-class queueing networks with
work-conserving service disciplines~\cite{dai1995_1,dai1995_2} have been derived
using these steps.
In many recent papers (see, for instance, \cite{gromoll2009}), very involved proofs have been considered
to demonstrate that the state of the stochastic network (under an appropriate space-time scaling)
 converges to a deterministic system whose evolution is represented by a differential equation of the form:
\begin{equation}\label{def:ode}
{d\over dt} x(t)= \delta(x(t)),
\end{equation}
where $\delta$ is the drift of the stochastic process and $x(t)$ corresponds to a fluid model.
It is worth mentioning that the proof of such a convergence found in \cite{gromoll2009} holds for state processes that are not necessarily Markovian.

It turns out however that such a convergence does not hold in general when the drift cannot be extended to a continuous function on $\mathbb R_+^N$. When the drift vector-field is discontinuous, the trajectories
of a fluid equivalent system enter sliding modes and the differential equation \label{def:ode} has to be replaced by a new dynamical system defined piece-wise by differential equations ${d\over dt} x(t)= \tilde \delta(x(t))$, where
$\tilde \delta$ is a convex combination of drifts of points in the neighborhood of the discontinuity. Such a phenomenon
was already emphasized in \cite{fayolle1995} where $\tilde \delta$ was called the "second vector-field" (see also \cite{meyn2007,fort2008}).
Unfortunately, $\tilde \delta$ is difficult to compute in many cases (it depends crucially on the statistical assumptions made) and has not been characterized in general.

\

Our contribution is three-fold. First, we present a way to construct a Lyapunov function
using the deterministic differential equations driving the fluid limits dynamics when the drift is continuous. This is done without proving that a deterministic differential equation indeed represents the behavior of a scaled version of the underlying process. Such proofs have been necessary so far while our approach is more direct and elementary for proving stochastic stability and gives another simple understanding of the meaning of such fluid limits for obtaining the stochastic stability. Moreover, it provides a systematical way of finding Lyapunov functions for stochastic systems. The advantage of finding a Lyapunov function explicitly is that it potentially gives much more precise information on the nature of the convergence of the process towards its stationary regime~\cite{meyn1993}. We also show that in the case of conservative drifts, the complexity of the problem can be considerably reduced.
We then give a counterpart result, leading to instability.

Second, we give general necessary conditions for stability in the case of discontinuous drifts.
These conditions have a natural geometric interpretation.

Third, we use these conditions to get a sharp geometric characterization of the stability set in case of piece-wise constant drifts in dimension~$2$. We give in particular an algorithm allowing to conclude if the process is stable or not, for all fixed birth rate parameters outside a set of dimension $1$.


\

The paper is organized as follows. In Section~\ref{sec:model} we describe the model in detail
and discuss the methodology used in the subsequent analysis.
In Section~\ref{sec:cont} we examine the case when the drift vector-field can be extended to a continuous function. Section~\ref{sec:discont} is devoted to deriving sufficient stability conditions
in the case of discontinuous drifts. We start by showing that fluid limits in this case are cumbersome and then proceed to presenting our approach in a generic scenario. In Section~\ref{sec:dim2}, we show how this approach may be applied to the processes in dimension~$2$ with piece-wise constant drifts in order to obtain a sharp geometric characterization of the stability region. Section \ref{sec:examples}
illustrates our various results and show that our sufficient conditions are not necessary in dimension 3.
Section \ref{sec:concl} concludes the paper.

\section{The model}\label{sec:model}


Let $N$ be an integer greater than $1$.
We denote by $\mathbb A^N_{+}$ the positive orthant of $\mathbb A^N$ (where $\mathbb A$ in this paper will be $\mathbb Z$ or $\mathbb R$) while
$\mathbb A^N_{+,*}$ stands for $\mathbb A^N_{+} \setminus \{ 0\}$.

Let $e_i$ be the vector of $\mathbb Z^N_{+}$ defined by
$(e_i)_i=1, \ (e_i)_j=0, j\neq i$.
If not specified otherwise, $|~|$ denotes the usual Euclidian norm.
The notation $x \le y$ is  used for the coordinate-wise ordering:
$\forall i, \ x_i \le y_i$, and we denote by $\langle x,y\rangle$ the usual scalar product of two vectors
in $\mathbb R^N$.
A process $X$ or a trajectory $u$ started in $x$ at time $0$ will be respectively denoted by $X^x$ and $u^x$.

\

Assume that $X$ is a continuous-time Markov process on $\mathbb Z_+^d$ with
the following transition rates:
$$q(x,x+e_i)=\lambda_i(x),$$
$$q(x,x-e_i)=\phi_i(x),$$
where ${\lambda} = (\lambda_i)_{i=1 \ldots d}$ and ${\phi}=(\phi_i)_{i=1\ldots d}$ are vectors consisting
of positive $0$-homogeneous functions from $\R^d_+$ to $\R$.

The drift function
$\delta=(\delta_1,\dots,\delta_N)=\lambda-\phi$ is bounded, which guarantees that the
process~$X$ is non-explosive. Hence we may assume that $X$ and all other stochastic
processes treated in the sequel have paths in the space $D = D(\R_+,\Z_+^N)$ of
right-continuous functions from $\R_+$ to $\Z_+^N$ with finite left limits.
Recall that a stochastic process with paths in~$D$ can be viewed as
a random element on the measurable space $(D, \D)$, where $\D$
denotes the Borel $\sigma$-algebra generated by the standard
Skorokhod topology~\cite{kallenberg2002}.

We are interested in conditions on the drift vector-field ${\delta}$ ensuring that the process is either stable (recurrent) or unstable (transient or null-recurrent). In subsequent sections we shall find such conditions with the use of the so-called Foster-Lyapunov criterion (see, e.g.~\cite{tweedie1978,fayolle1995,meyn1993}).

\section{Smooth drift} \label{sec:cont}

\subsection{Stability conditions}

Let $(ODE)_x$ be the following deterministic differential equation
\begin{equation}
{d \over dt} u(t)=\delta(u(t)).\label{eq:odex}
\end{equation}
$$u(0)=x.$$

We denote by $u^x$ a solution with an initial condition $x$.
Define~$\S=\{ x : |x| =1\}$.

\begin{theorem}

Assume that $\delta$ is a continuously differentiable function from $\mathbb R_+^N$ to $\mathbb R_+^N$.
Assume that for all $x \in \S$, there exists a solution of $(ODE)_x$ such that $x(t) =0$
for all $t \ge T_x$ where $T_x < \infty$. Assume in addition that
\begin{equation} \label{eq:unif_cont}
\sup_{x \in \S} T_x < \infty.
\end{equation}
Then $X$ is positive recurrent and $x \mapsto T_x$ is a Lyapunov function.
\end{theorem}

\bp
Note first that the homogeneity of the drift implies that if $u^x(t)$ is a solution of $(ODE)_x$,
then we can define $u^{Kx}(t)=Ku^x(t/K)$ as a solution of $(ODE)_{Kx}$. Indeed:
$$u^{Kx}(0)=Kx,$$
$$ {d \over dt} u^{Kx}(t)=K {d \over dt}u^{x}(t/K) =\delta(u^x(t/K))= \delta(Ku^x(t/K))=\delta (u^{Kx}(t)).$$

Because $\delta$ is $C^1$ on $\mathbb R^N_{*}$, the flow $(t,x)\to u^x(t)$ is $C^1$ on $\mathbb R \times \mathbb R^N_*$, i.e. continuously differentiable in $t$ and $x$ (see for instance Theorem 1, page 299 in \cite{hirsch1974}).

Define by $F(x) = T_x$ the time needed for $u^x$ to hit $0$. We are going to show that $F$ is a suitable Lyapunov function for proving positive recurrence of $X$. Due to the assumptions of the Theorem, $F$ is a positive finite function. It is easy to see that
$${d \over dt} F(u^{x}(t)) = -1 \quad \text{for all} \quad t < T_x.$$
Indeed, the difference between $F(u^x(t+h))$ and $F(u^x(t))$ is negative and is equal in absolute value to the time needed to reach $u^x(t+h)$ from $u^x(t)$, which is exactly $h$. Hence, the latter equality follows.

Moreover since $u^{Kx}(t)=Ku^x(t/K)$, it follows that $F(Kx)=KF(x)$, i.e. $F$ is a $1$-homogeneous function. This, in particular, implies that $F(x) \to \infty$, when $|x| \to \infty$.
The drift of F is given by
$$\Delta F({x})= \sum_{{y}} q(x,y)(F({y})-F({x})). $$

It will be convenient to approximate the variations of $F$ using its gradient.
We thus need to show that $F$ is differentiable. For that purpose, define for each $\kappa > 0$ the functions
$$
T_{x, \kappa} = \inf\{t: |u^x(t)| \le \kappa\}.
$$
It is clear that $T_{x, \kappa} \to T_x$ for each fixed $x$ as $\kappa \to 0$. It is also clear that for each fixed $x$ functions $T_{x, \kappa}$ increase when $\kappa$ decreases. Note that, due to the continuity of $u^x(t)$ in $t$,
$$
|u^x(T_{x, \kappa})| = \kappa.
$$
Examine this equality: $u^x(t)$ is a differentiable function, while  the norm of a differentiable function is also differentiable. Hence, we conclude that $T_{x, \kappa}$ is differentiable for all values of $\kappa > 0$. In order to prove the differentiability of $T_x$ (and $F(x)$), it remains to show that the convergence $T_{x, \kappa} \to T_x$ is uniform in $x$. It follows from the following sequence of equalities:
\begin{eqnarray*}
\sup_x |T_{x, \kappa} - T_x| = \sup_{x: |x| = \kappa} |T_x| = \kappa \cdotp \sup_{x \in \S} |T_x| \to 0
\end{eqnarray*}
as $\kappa \to 0$, due to (\ref{eq:unif_cont}). The last equality in the sequence follows from $1$-homogeneity of $T_x$.

As $F$ is $1$-homogeneous and differentiable, we have
\begin{equation} \label{eq:approx}
\Delta F({x}) \le \left\langle\nabla F({x}),\delta(x) \right\rangle + \varepsilon,
\end{equation}
for $|{x}|$ large enough, where $\varepsilon > 0$ may be made
sufficiently small.

Observing that
 $$\dfrac{dF(x(t))}{dt} = \left\langle \nabla F({x}(t)), \dfrac{d{x}(t)}{dt} \right\rangle = \left\langle\nabla
F({x(t)}), \delta(x(t))\right\rangle,$$
we obtain, taking $t \to 0$, that
$$\left\langle \nabla F(x),\delta(x) \right\rangle \le -\epsilon,$$
For $|x|$ large enough.
Summarizing, the function $F$ is such that $F({x}) \to \infty$ as $|{x}| \to \infty$ and
that $\Delta F(x) \le -\epsilon$,  for $|x| $ large enough.
After an examination of the Foster-Lyapunov criterion, it can be concluded that $X$ is positive recurrent.

\ep

\begin{remark}
Note that our approach is very similar to that of the fluid-limits approximation as it is developed for instance in \cite{fort2008}. However, we do not need to prove that the differential equation we look at is indeed the one that represents the behavior of our birth-and-death process on the fluid scale. Instead, we have clearly used that the differential equation $(ODE)_x$ is the characteristic equation of the PDE (partial differential equation) $\left\langle \nabla F(x),\delta(x) \right\rangle = -\epsilon$ which directly leads us to finding a Lyapunov function.
\end{remark}

\begin{remark}
 It is worth mentioning that our analysis is valid thanks to the $0$-homogeneity assumption.
 Without this assumption, the stability of the birth and death processes {\bf cannot} be described using the differential equation (\ref{eq:odex}). A counter-example (with continuous drifts) can be found
 in \cite{borst2008}.
\end{remark}

\begin{remark}
We believe that our construction could be considered for continuous drifts vector-fields using approximation arguments.
This falls however out of the scope of this paper.
\end{remark}

\subsection{Instability conditions}
In this section, we consider a reverse statement establishing instability
relying on the previously considered dynamical system. 
It is much more challenging to state generic instability conditions based on the ODE~(\ref{eq:odex}),
without a direct use of fluid limits i.e. without needing to prove the convergence of a scaled version of the process towards the trajectories of the ODE (see \cite{meyn1995} for the construction of a Lyapunov
function proving the transience of multi-class queuing networks with routing).

We use hereafter explicitly the convergence to the fluid limits and refer to \cite{gromoll2009} (see also in the discrete time setting \cite{fort2008})
for a proof of this result. The next theorem 
is hence essentially a combination of proving the convergence of the scaled process
 and the extended version of the (instability part of) Foster-Lyapunov criterion (see for instance Theorem 2.2.7 in \cite{fayolle1995}).

%


\begin{theorem}\label{theo:inst}
Assume that $\delta$ is a Lipschitz function outside a neighborhood of the origin. Assume further that there exists a strictly positive time $T$ and a number $a>1$ such that
for all $x$, with $|x|=1$, a solution $u^x$ of $(ODE)_x$ is defined on an interval $[0,\tau_x]$
with $\tau_x > T$
and verifies $|u^x(T)|\ge a$.
Finally suppose that for all times $t$ and points $x$
$$\lim_{K \to \infty} P(\sup_x \left|{X^{Kx}(Kt))\over Kt}\right| > \epsilon )=0.$$
Then $X$ is transient.
\end{theorem}

\bp
Note first that the Lipschitz condition ensures that the ODE~(\ref{eq:odex}) has a unique solution $u^x$ for each $x \in \mathbb R_{+,*}^N$ on an interval $[0,\tau_x]$. Furthermore, since $0$ is necessarily a stable point if it is an equilibrium point, the conditions of the theorem
imply that the trajectories $u^x(t)$ did not hit $0$ before time $T$.
Using \cite{gromoll2009}[Theorem 4.1] (which actually needs only that $\delta$ is continuous), we have the convergence of the process towards its fluid limit in the sense that
$$\lim_{K \to \infty} P\left(\left|{X^{Kx}(Kt) \over K}-u^x(t)\right| \ge \epsilon \right) =1,$$
for each interval $[0,t]$ included in the interval $[0,\tau_x]$ where the ODE has a solution.

This implies that for all $x$, there exists $K_x$ such that
$$E|X^{K_x x}(K_x T)| - K_x |x| \ge (u^x(T)-1) K_x +\epsilon >0.$$
Hence
$$E|X^{K_x x}(K_x T)| - K_x |x| \ge (a-1) K_x.$$
We now need to prove that $\sup K_x =K <\infty$ to be able to make use of the extended Foster--Lyapunov criterion (see for instance \cite{fayolle1995}[Theorem 2.2.7]).
Denote by $\B_{\eta}(x)$ a ball with radius $\eta$ and center $x$.

Fix $\eta>0$.
Using the martingale decomposition, $X^{Kx}$ can be decomposed as
\begin{eqnarray*}
\frac{X^{Kx}(t)}{K}&=& {K x \over K} +  {1 \over K} \int_{0}^{Kt} \delta(X^{Kx}(s))ds +  { M_{Kt} \over K}\\
&= &{x} +  {1 \over K} \int_{0}^{Kt} \delta(X^{Kx}(s))ds  +  { M_{Kt} \over K},
\end{eqnarray*}
where $M$ is a martingale that satisfying:
$$E\left(\sup_{0\le s\le t} \frac{M_{Kt}}{K}\right) \le A \left({t \over K} \right)^{1/2} \le \eta .$$
for $K$ large enough and $A$ being a positive constant.
Defining now $f_K(t)=E[\sup_{0 \le s \le T} {1 \over K} |X^{Kx}(Ks) -X^{Ky}(Ks)|]$, we get that:
\begin{eqnarray*}
f_K(t) &=&   |x-y|+  2\eta + {1\over K} E\left[\sup_{0\le s \le T}  \left| \int_{0}^{Ks} \delta(X^{Kx}(u))du
 -\int_{0}^{Ks} \delta(X^{Kx}(u))du\right| \right], \\
&\le& 3 \eta  + E\left[\sup_{0\le s \le T}\left|\int_{0}^{s} (\delta({X^{Kx}(Ku) \over K})- \delta({X^{Ky}(Ku) \over K} )du\right|\right],
 \end{eqnarray*}
where we have successively used the convergence of the martingale part of the process to $0$ and the homogeneity of the drift.
We can now condition on the event $\A_\epsilon=\{\forall s; X^{Kx}(s) \notin \B_\epsilon(0)\} \cap \{\forall s : X^{Ky}(s) \notin \B_\epsilon(0)\}$ and use the Lipschitz condition. There exist two constants $B$ and $C$ such that:
\begin{eqnarray*}
f_K(t)  &\le& 3 \eta + B E \int_{0}^T\left| {X^{Kx}(Ku) \over K}- {X^{Ky}(Ku) \over K} \right| du + C P(\A_\epsilon).
 \end{eqnarray*}
Using the theorem assumptions, we get that there exists $K_0$ such that for all $K \ge K_0$, $P(\A_\epsilon) \le \eta$.
We now conclude from Gronwall's lemma that there exists a constant $D$ such that for all $y \in \B_{\eta}(x)$:
$$f_K(t) \le D \eta ,$$
which in turn implies that $\sup K_x =K <\infty$.
The last assertion concludes the proof.



%


\ep

\subsection{Gradient systems}

Finding an explicit form for the function $F$ might be difficult in general. However, under slightly stronger assumptions on the drift function $\delta$, we can construct an explicit Lyapunov function directly from the vector-field $\delta$. This is a well known fact in the theory of deterministic dynamical systems.

\begin{proposition} \label{prop:potential}
Assume that $\delta$ is a conservative vector-field i.e.,  $\delta= - \nabla V$,
and assume that
\begin{equation} \label{cond:V}
V(x) \ge a>0.
\end{equation}
Assume further that there exists $\epsilon >0$, such that
\begin{equation} \label{cond:d}
|\delta(x)| \ge \epsilon \quad \text{for all} \quad x.
\end{equation}
Then $V(x)$ is a Lyapunov function and $X$ is positive recurrent.
\end{proposition}

\bp
Note first that since $\delta$ is $0$-homogeneous, $V$ is $1$-homogeneous.
Using the $1$-homogeneity, there exists $\kappa < \epsilon^2$ such that for $|x|$ large enough, we can estimate the drift of $V$ by:
\begin{eqnarray*}
\Delta V(x)&\le& \langle \nabla  V(x), \delta(x)  \rangle + \kappa\\
&=&  \langle \nabla ~ V(x), - \nabla V(x) \rangle + \kappa\\
&= & -|\nabla ~ V(x)|^2  + \kappa \le -\epsilon'.
\end{eqnarray*}
Furthermore, $|V(x)| \to \infty$ for $|x| \to \infty$ since it is a $1$-homogeneous  and strictly positive
function.
We can therefore apply the Foster--Lyapunov criterion.
\ep

\begin{remark}
A vector-field $\delta(x_1,x_2)=(\delta_1(x_1,x_2), \delta_2(x_1,x_2))$ (on a completely connected set) is conservative if and only if ${d\over dx_1}\delta_2(x)={d\over dx_2}\delta_1(x).$
\end{remark}




 \

\section{Discontinuous drifts} \label{sec:discont}

\subsection{Complexity of the fluid limits}

So far we restricted ourselves to the case when the drift vector-field is continuous. The situation changes dramatically when this condition is dropped.
When the drift vector-field is discontinuous, the trajectories
of a fluid equivalent system near a point of discontinuity may enter "sliding modes" and the differential equation~(\ref{def:ode}) has to be replaced by a new dynamical system defined piece-wise by differential equations ${d\over dt} x(t)= \tilde \delta(x(t))$, where
$\tilde \delta$ is a convex combination of drifts around neighborhoods of the discontinuities. 

Let us give a simple example of this phenomenon.
Consider the following transitions with fixed birth rates  $\lambda_1,\lambda_2$ and death rates given by the following bandwidth allocation:
\begin{eqnarray} \label{eq:defcoupled}
\phi_1(x)=1_{x_2=0} + a_1 1_{x_2>0},
\phi_2(x)=1_{x_1=0} + a_2 1_{x_1>0}.
\end{eqnarray}

Suppose $\lambda_1 < a_1$ and $\lambda_2 < (1-\rho_1)+ a_2 \rho_1$. This condition is known to be sufficient for stability of such a model in dimension~$2$, and it has been obtained through different methods (see, for instance, \cite{cohen1983,fayolle1995,borst2008}).

The following proposition proved in \cite{robert2003} characterizes the fluid limit of $X$.

\begin{proposition}
The process ${X^K(Kt) \over K}$ converges in distribution when $K \to \infty$ towards a process
$x(t)$ satisfying the differential equations:
\begin{eqnarray}
{d\over dt} x(t)&=&\delta(x(t)), \mbox{ for } x(t)>0,\\
{d\over dt} x_1(t)&=&0, \mbox{ for } x_1(t)=0,\\
{d\over dt} x_2(t)&=& \lambda_2 - (1-\rho_1)+ a_2 \rho_1, \mbox{ for } x_1(t)=0.
\end{eqnarray}

\end{proposition}

The stability condition is easily interpreted when considering the convergence to $0$ of the obtained fluid limit.
This example shows however that even in a very simple case the fluid limit satisfies an equation different from (\ref{def:ode}).

In the next subsection we develop an approach that allows us to find stability conditions in the case of discontinuous drift vector-fields without the use of fluid approximation.

\


\subsection{Sufficient stability conditions} \label{sec:dv}

This section is devoted to identifying rather general conditions on the drift vector-field
ensuring stability even in the presence of discontinuities for $\delta$.
These conditions lead to useful geometric stability conditions in dimension 2, which are discussed in Section~\ref{sec:dim2}.

We start by considering a general vector-field of $0$-homogeneous drifts such that the number of discontinuities is finite. We construct a Lyapunov function by pasting together local Lyapunov functions and using a smoothing technique. This method was first used in \cite{dupuis1994}.

We need to introduce a few more notations.
Define by $\B_\varepsilon(x)$ a closed sphere with radius $\epsilon$ and center $x$.
For a point $x$, denote by $\D_{\varepsilon}(x)$ the set of drifts in a neighborhood of $x$, i.e.,
$$\D_{\varepsilon}(x)=\{ \delta(y): y \in B_\varepsilon(x) \}.$$
We then define by $\D^*_{a,\varepsilon}(x)$ the set of vectors
 $$\D^*_{a,\varepsilon}(x)=\left\{ \eta \in \mathbb R^d ~ : \langle \eta, v \rangle < - a, \forall v \in \D_{\varepsilon}(x) \cup\{-x\} \right\} .$$


We now state an assumption on the vector-field $\delta(x)$ that we shall prove to be sufficient to characterize the stability region of the process.


Assumption $(A_1)$: 
For all $x \neq 0$, there exists $\epsilon>0$ and $a >0$ such that:
 $$ \D^*_{a,\epsilon}(x) \neq \emptyset .$$.

The theorem below is the main result of this section.

\begin{theorem} \label{thm:discont}
Assumption $(A_1)$ implies that $X$ is positive recurrent.
\end{theorem}

Before presenting a rigorous proof, we would like to explain the result intuitively.
If the sets $\D_{\varepsilon}(x)$ are finite for all $x$, assumption (A1) may be better understood using a simple geometric interpretation.
Using Farkas' lemma, we can state that either assumption $(A_1)$ is true or $x$
is in the cone induced by the vectors of $\D_{\varepsilon}(x)$, i.e. there exist non-negative weights
$\alpha_i$ such that
$$\sum_{i \in \I} \alpha_i \delta^i =x.$$
It is hence natural to expect that if $x$ is never contained in the cone induced by the drifts $\delta(y)$
at points $y$ close to $x$ (which is exactly assumption (A1)), then the process is stable.

\bp

Consider first the function $H$ and the vector-field $v_u$ constructed by
$$v_u(x)=\arg\max_{\eta \in \D^{*}_{u,a}(x)} \langle x,\eta \rangle .$$
$$H(x,u)=\max_{\eta \in \D^{*}_{u,a}(x)} \langle x,\eta \rangle .$$
The function $H$ is a natural candidate for a Lyapunov function but it is discontinuous which complicates
drastically the drift calculations and precludes having a negative drift in all points.
We overcome this difficulty by considering a smoothed version of $H$.
Let $\kappa_\epsilon$ be a $C^{\infty}$-probability density supported on the sphere $\B_\varepsilon(0)$ and introduce
$$F(x)=\int_{u \in \B(\epsilon)} H(x,u)\kappa_{\epsilon}(u) du.$$
Then $F$ is clearly a $C^{\infty}$-function. We shall now prove that $F$ is a suitable Lyapunov function.
First, notice that by assumption $(A_1)$, $H(x,u) \ge a$, and hence, $F(x) \ge a$ for all $x$. This, together with the observation that the function $F$ is $1$-homogeneous, implies that $F(x) \to \infty$ as $|x| \to \infty$.

Note further that the compactness of the sphere guarantees the existence of $u \le \epsilon$ (arbitrarily small) such that for all
$x$,
\begin{equation}\label{ineq:A}
|\nabla F(x)- v_u(x)| \le -\epsilon.
\end{equation}

Using assumption $(A_1)$ again, we get that
\begin{equation}\label{ineq:B}
\langle \delta(x), v_u(x) \rangle \le -a.
 \end{equation}

Hence using \ref{ineq:A},
\begin{equation}\label{ineq:C}
 \langle \delta(x), \nabla F(x) \rangle \le \langle \delta(x), v_u(x) \rangle + C \epsilon.
\end{equation}

Again using the fact that $F$ is $1$-homogeneous, and combining \ref{ineq:B} and \ref{ineq:C} we get that
$$\Delta F(x) =\sum_y q(x,y)(F(y)-F(x)) \le \langle \delta(x), \nabla F(x) \rangle + \epsilon',$$ for large $|x|$.
Hence, $F$ is a Lyapunov function.

\ep

\section{Piece-wise constant drift in dimension 2} \label{sec:dim2}

In this subsection we apply the general result of Theorem \ref{thm:discont} to a particular case of a discontinuous drift function. We assume that the state space of the underlying process is $\mathbb N^2$ and that the rate functions are piece-wise constant. Together with the assumption that the rate functions are $0$-homogeneous, this means that there is a finite number of cones where rate functions are constant.

We start be introducing some notations that will be used throughout this section. Assume that there are $N$ vectors $v_1, ..., v_N$ such that $v_1 = e_1, v_2 = e_2$ (where $e_1$ and $e_2$ are vectors co-directed with one of the axes) and such that $\delta(x) = \delta^k$ for any $x = Av_k + Bv_{k+1}$. This means that the drift at any point of the cone defined by $v_k$ and $v_{k+1}$ is equal to $\delta^k$.

Note that we do not require the vectors $v_k$ to be different. (This means that a cone reduces to a line when two consecutive vectors $v_k$ and $v_{k+1}$ are equal.)

Introduce also certain sets that will be crucial for the definition of the stability region. For each $k=1,..., N-1$, let
$$
U_k^1 = \left\{\delta : \delta^k=  Av_k + Bv_{k+1} \quad \text{for some} \quad A \ge 0, B \ge 0, A+B > 0\right\}
$$
and for each $k=2,..., N$, let
$$
U_k^2 = \left\{\delta:   \alpha \delta^k + (1-\alpha) \delta^{k-1}= A v_k   \quad \text{for some} \quad A \ge 0, \alpha \in [0,1] \right\}.
$$

We are ready to state the main result of this section.

\begin{theorem} \label{thm:2dim}
Assume that
$$
\delta \in \SS= \bigcap\limits_{k=1}^{N-1} \left(\overline{U_k^1} \bigcap \overline{U_k^2}\right).
$$
Then the Markov process $X$ is positive recurrent.
Conversely, if $\delta$ belongs to the interior of the complement of $\SS$, then $X$ is transient or null-recurrent.
\end{theorem}

\bp

We start by proving the first part of the theorem. In order to do this, let us verify the conditions of Theorem~\ref{thm:discont}. It is clear that only two situations are possible:

(i) the vector $x$ belongs to the interior of a cone defined by vectors $v_k$ and $v_{k+1}$. In this case $x = Av_k + Bv_{k+1}$ for some $A > 0, B >0$.

(ii)the vector $x$ is collinear to a vector $v_k$. In this case $x = A v_k$ for some $A>0$.

Consider these two situations separately. In case $(i)$, thanks to Theorem~\ref{thm:discont}, we need to show the existence of a vector $\eta$ such that
$$
\langle \eta, Av_k + Bv_{k+1} \rangle > 0
$$
and
$$
\langle \eta, \delta^k \rangle < 0.
$$
In geometric terms this means that there exists a vector $\eta$ such that vectors $Av_k + Bv_{k+1}$ and $\delta^k$ belong to different half-planes separated by the line normal to vector $\eta$. It is easy to see that the existence of such a vector is guaranteed by the fact that $\delta \notin U_k^1$.

Consider now situation $(ii)$. Applying again Theorem~\ref{thm:discont}, we see that one needs to show the existence of a vector $\eta$ such that
$$
\langle \eta, Av_k \rangle > 0,
$$
$$
\langle \eta, \delta^k \rangle < 0
$$
and
$$
\langle \eta, \delta^{k-1} \rangle < 0.
$$
If we interpret this again in geometric terms, it is equivalent to the existence of a vector $\eta$ such that vectors $-Av_k$, $\delta^k$ and $\delta^{k-1}$ belong to the same half-plane defined by the line normal to $\eta$. Direct computations show that it follows from the fact that $\delta \notin U_k^2$.
The proof of the positive recurrence under the assumption that $\delta \in \SS$ is now complete.

\

Let us now show that if $
\delta \in int \left( \bigcup_{k=1}^{N-1} U_k^1 \bigcup U_k^2 \right)$,
then the Markov process $X$ is not positive recurrent.
We are going to prove that under the given conditions,
the process is actually not rate stable which prevents stability.
Assume that the process is started in a cone $k$ and that $\delta \in U_k^1$.
The strong law of large numbers (SLLN) then implies that with a positive probability the process stays in the cone $\{ A v_{k} +B v_{k+1} , A \ge 0, B \ge 0\}$.

Assume that $X$ is positive recurrent. This implies that ${X^x_t \over t} \to 0$, almost surely.
Using the Martingale decomposition of the process,
$${X^x_t \over t}= {x \over t} + {M_t \over t} + {1 \over t}\int_{0}^t \delta(X^x_s)ds.$$
Using the boundedness of the transitions, the martingale $M_t$ is such that $E[M_t^2] \le C t$,
which implies the convergence in $L^2$ and in probability of $M_t \over t$ to $0$, which in turns implies the almost sure convergence along a subsequence. 
Conditioning on the fact that the process stays in the cone we obtain ${1 \over t_n} \int_0^{t_n} \delta(X^x(s))ds \to 0, t_n \to \infty$, which combined with the ergodic theorem for the positive recurrent Markov process $X$ implies that there exists $\alpha \ge 0$ such that
$$0 = \alpha \delta^{k-1}+ (1-\alpha) \delta^k,$$
which contradicts $\delta \in U^1_k$.

Suppose now that $\delta \in U^2_k$. In this case the SLLN implies that with a strictly positive probability, the process stays in the set $\{ A v_{k-1} +B v_{k} , A \ge 0, B \ge 0\}\cup \{ A v_{k} +B v_{k+1} , A \ge 0, B \ge 0\}$.
Proceeding similarly as in the previous case, there exists $\tilde \alpha$ such that:
$$0 = \tilde \alpha \delta^{k}+ (1-\tilde \alpha) \delta^{k+1},$$
which contradicts $\delta \in U^2_k$.

\ep

\subsection{Fluid limits}

A very interesting situation occurs within the framework of this section when there exist $k$ and $\alpha \in (0,1)$ such that $\alpha \delta^k + (1 - \alpha) \delta^{k+1} = A v_k$ for some $A \neq 0$.
In this case we know that the fluid limits\footnote{for the existence of fluid limits, we refer to \cite{robert2003} and to \cite{fort2008} Theorem 1.2, in discret time} with an initial state in the cone defined by $v_{k-1}$ and $v_k$ or in the cone defined by $v_k$ and $v_{k+1}$ enter a so-called "sliding mode" with their trajectory reaching the ray defined by $v_k$ after a finite time and not leaving it after this time.
The new drift $\tilde \delta$ obtained during the sliding mode on $v_k$ must be a convex combination of $\delta^k$ and $\delta^{k+1}$ and must also be collinear with $v_k$. Hence, we can explicitly calculate $\tilde \delta$ by solving the following system in $\alpha$ and $A$:
\begin{equation}\label{eq:tdrift}
\alpha \delta^k + (1 - \alpha) \delta^{k+1} = A v_k.
\end{equation}
The existence of a solution with a strictly negative $A$ is necessary to get a stable system while the existence of a solution
with strictly positive $A$ is sufficient to get instability of the process, which corresponds respectively to the case where
the fluid limits converge to zero or infinity.

This is a very particular scenario, as in the case of a dimension higher than $2$, equation (\ref{eq:tdrift})
is generally underdetermined and stability conditions cannot be characterized directly.

\begin{figure}[!]
\includegraphics[scale=0.25]{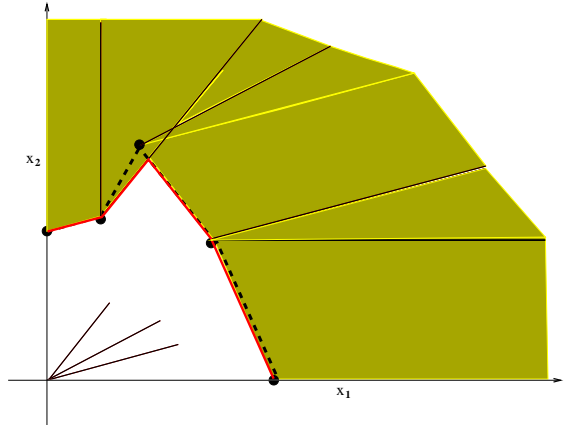}
\caption{Generic construction of the stability set for piece-wise constant drifts in dimension $2$.}
\label{pic:ex2}
\end{figure}

\subsection{Algorithm.} \label{subsec:alg}

Recall that we defined the drift vector $\delta$ to be equal to $\lambda(x) - \phi(x)$. It is often the case in queueing and telecommunications applications that the function $\lambda$ (representing the arrival rate) is assumed to be a constant and the question is for which values of $\lambda$ the system under consideration is going to be stable. We present here an algorithm to construct the stability set when $\lambda$ is fixed. Define (within the description of the algorithm and the examples later on) $\delta^k = \lambda - \psi_k$. The algorithm is given below.

\begin{itemize}
\item
step 1: Draw the points representing various values of $\psi_k$.

\item
step 2 : Connect $\psi_1$ to $\psi_2$, $\psi_2$ to $\psi_3$, etc.

\item
step 3: For each $k$, draw the cone defined by vectors $v_k$ and $v_{k+1}$ based on point $\psi_k$.

The compact set obtained is the stability region.

\end{itemize}

A generic illustration of the previous algorithm is given in Figure \ref{pic:ex2}.

\section{Examples} \label{sec:examples}

This section contains a number of examples illustrating the use of our results in various settings.

\subsection{Continuous drifts}

\begin{example}
Consider a bandwidth sharing network representing a data communication network as described in \cite{bonald2006} with the following constraints on routing and capacity:
$$\lambda \in \A,$$
$$\phi \in \C,$$
where $\A$ and $\C$ are two convex sets.

The following type of policies has been considered for many models in performance analysis (see for instance\cite{meyn2007}).
Assume that the vector of traffic intensities and the bandwidth allocation are chosen such that
$$\delta(x)=\arg\max_{{\lambda,\eta} \in \A \times \C } \langle x,\lambda -\eta \rangle=-\nabla \delta^*_{\A - \C},$$
where $\delta^*_S$ is the Minkowski function associated with a convex set $\S$, i.e.,
$$\delta^*_S(x)=\max_{u \in \S}  \langle x,u \rangle,$$
and is also the Fenchel-Legendre transform of the support function $1_{u \in \S}$.
It is now a matter of routine arguments to show that conditions \ref{cond:V} and \ref{cond:d} are satisfied if the interior of $ \A \cap \C$ is not the empty set.
\end{example}

\begin{figure}[!]
\includegraphics[scale=0.3]{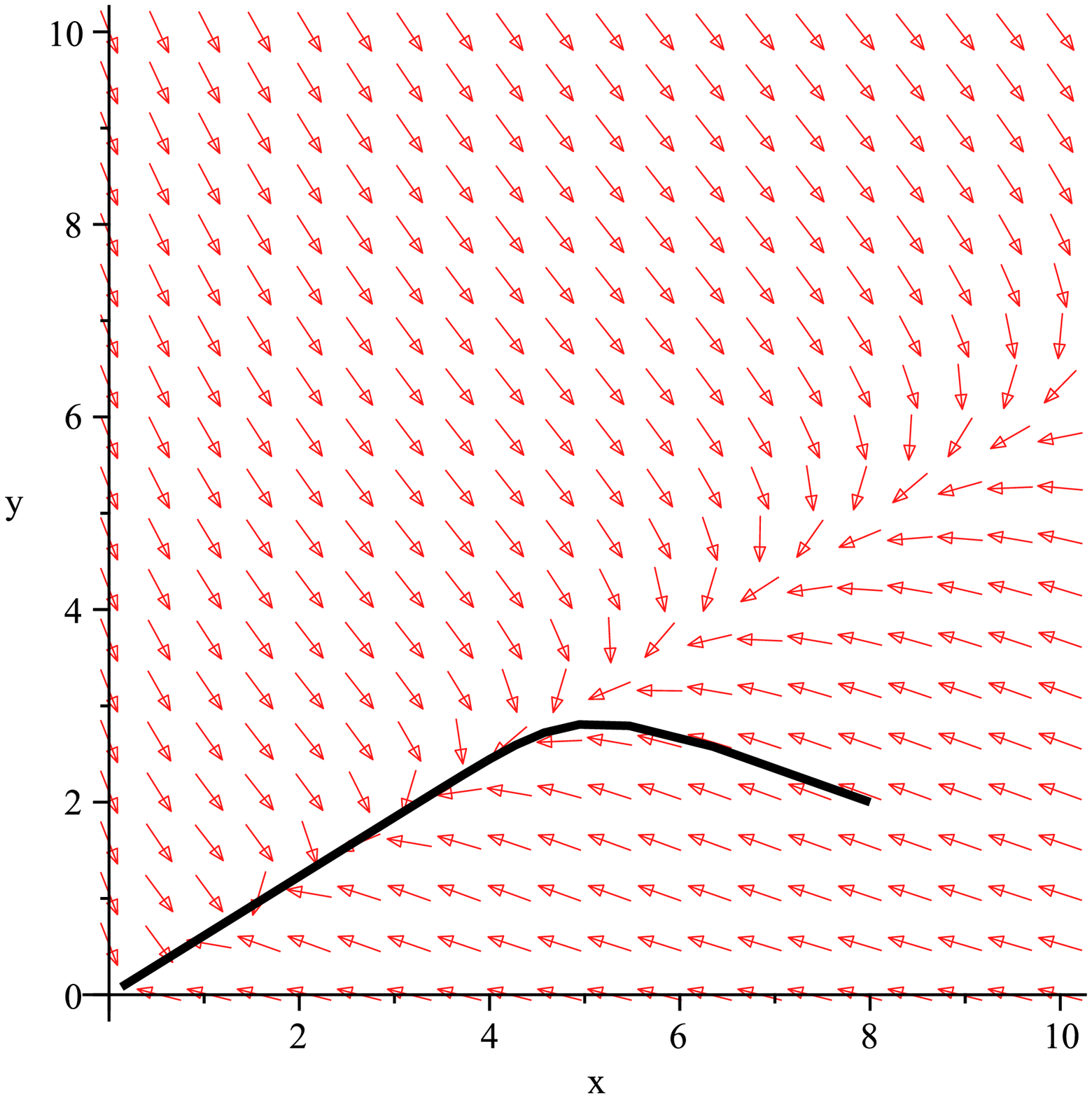}
\includegraphics[scale=0.3]{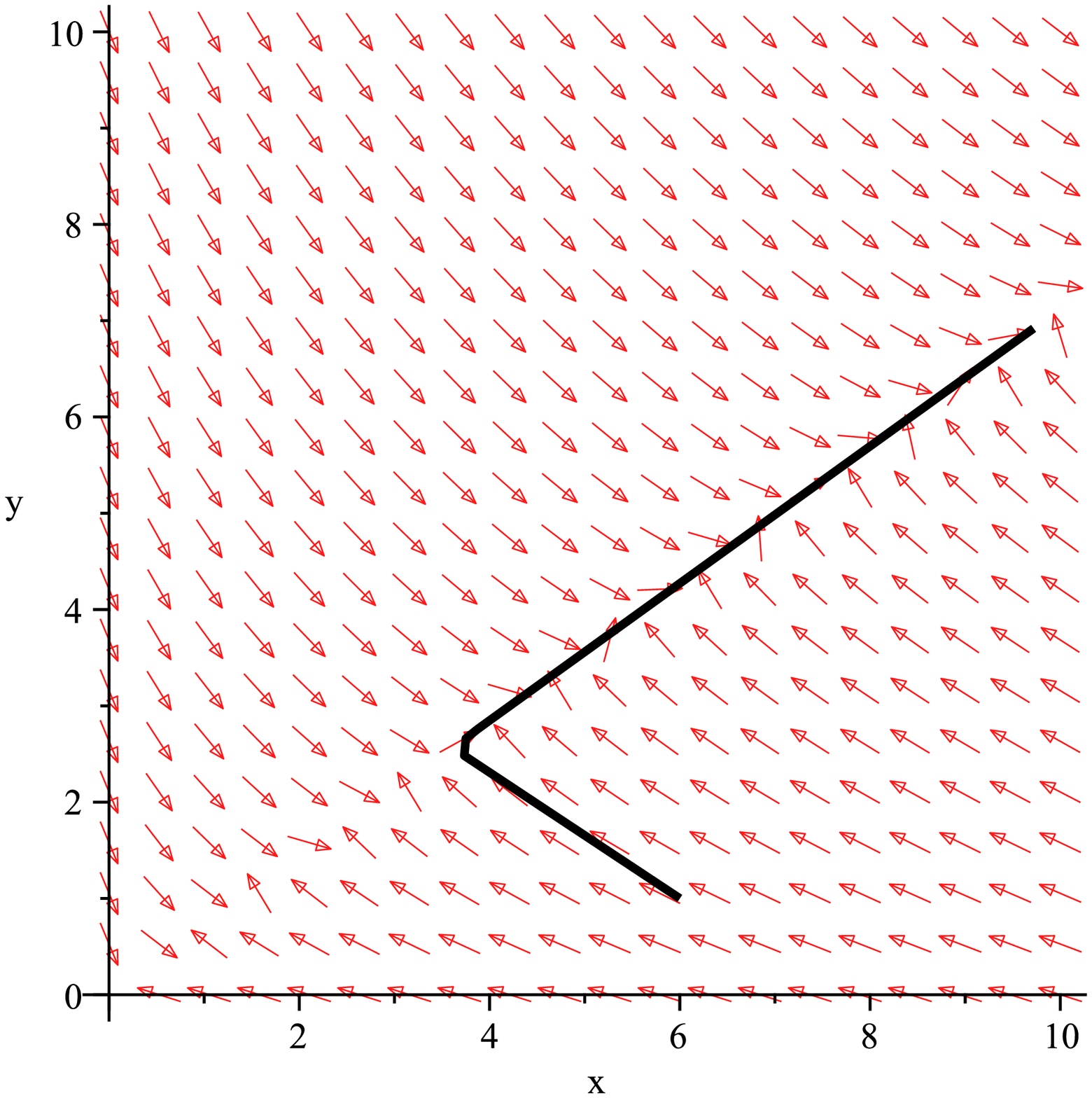}
\caption{Example 3 : Drift vector fields for birth rates $(0.4, 0.8)$ and $(0.5, 0.8)$.}
\label{fig:1}
\end{figure}

The results obtained in Section~\ref{sec:cont} allow to study numerically the positive recurrence of processes
 even when the drifts are too complicated to get an explicit solution for the associated ODE.

\begin{example}
An example of a wireless network with two types of users competing for the same bandwidth could
lead to the following deaths rates (using Shannon's formula and a state-dependent allocation policy):
$$\phi_1(x)=log\left(1+{x_1/|x|\over N + x_2/|x|}\right),$$
$$\phi_2(x)=log\left(1+{x_2/|x|\over N + x_1/|x|}\right),$$
where $N$ is the thermal noise.

Let us consider two possible vectors of arrival (birth) rate $(0.4, 0.8)$ and $(0.5, 0.8)$.
The associated (ODE) can be solved numerically for any value $\lambda_1, \lambda_2$
allowing to conclude for the positive recurrence of the process in the first case and the transience in the second case, as shown in figure \ref{fig:1} using the following properties of the trajectories of the
ODE:
\begin{itemize}
\item
In the first case, all trajectories started from any state on the sphere hit $0$ in a bounded time,

\item
In the second case, all trajectories started from the sphere do no reach a sufficiently small neighborhood of $0$, from which we can conclude that all fluid-limits solutions stay outside of a ball of radius $\epsilon$ and center $0$. Moreover all trajectories do reach a state of norm bigger than $1$ before a finite time $T$.

\end{itemize}

\end{example}

\subsection{Discontinuous drifts}

We will give here a few examples where Theorem~\ref{thm:2dim} provides interesting results. 

\begin{figure}[!]
\includegraphics[scale=0.5]{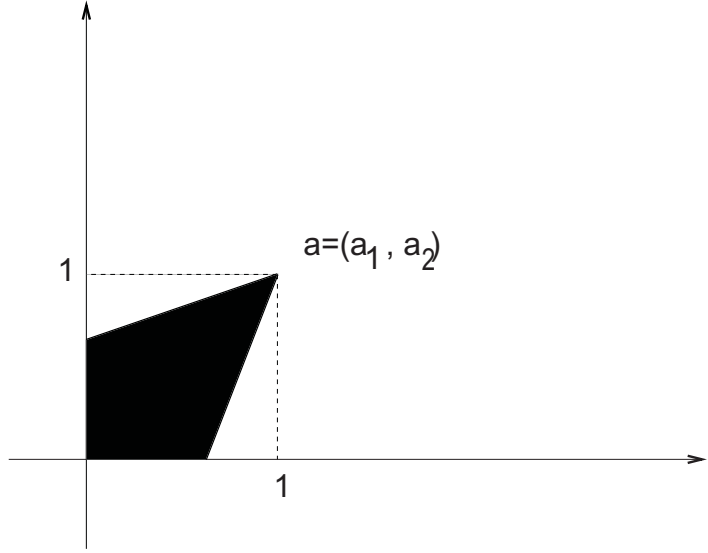}
\caption{Stability region (set of birth rates $\lambda$) for the 2 coupled processors of example 4}
\label{pic:ex3}
\end{figure}

\begin{example}

We describe here how Theorem~\ref{thm:2dim} can be used to obtain the well-known stability results for the so-called coupled-processors problem.
Consider  the allocation described as the most basic example with discontinuous drifts in Section
 \ref{sec:discont}:
  \begin{eqnarray*}
\phi_1(x)=1_{x_2=0} + a_1 1_{x_2>0},
\phi_2(x)=1_{x_1=0} + a_2 1_{x_1>0}.
\end{eqnarray*}
It is clear that in this case the algorithm of Subsection~\ref{subsec:alg} allows us to recover the well-known stability region for this problem \cite{cohen1983,fayolle1995} (see Figure~\ref{pic:ex3}):
$\lambda_1 < a_1$ and $\lambda_2 < (1-\rho_1)+ a_2 \rho_1$ or $\lambda_2 < a_2$ and $\lambda_1 < (1-\rho_2)+ a_1 \rho_2$, with $\rho_i={\lambda_i \over \mu_i}$.
\end{example}

\begin{remark}
Based on the previous example, one may think at first sight that the stability region is the smallest convex set containing all the vectors $\phi(x)$ for $x$ describing the state space. However, this is clearly not the case as it is illustrated in the generic example of Figure \ref{pic:ex2}. In fact, the stability set may not be convex.
\end{remark}

\begin{example}
We now look at the same type of death rates as in the previous example but with different birth rates.
Consider for instance a queuing or manufacturing system with $2$ processors and two types of traffic:
\begin{itemize}
\item
some dedicated traffic arriving with intensity $\lambda_i$ to processor $i$,

\item
some flexible traffic with intensity $\nu$ that can be routed to either processor $1$ or $2$ depending on the congestion level of both processors.
\end{itemize}
Assume that the flexible traffic is actually routed to the processor with the smallest number of jobs in process
(and to processor $1$, say, if the processors are equally loaded, this last assumption having no impact on the stability conditions).
 Assume further that the arrivals of jobs of each type of traffic are following a Poisson process (independent of each other and of everything else) and that the processing times are exponentially distributed. The presence of interference or switching costs between
different type of tasks raise the allocation of service (or death rates) described in the previous example with $a_i <1$.

Using the notations of the previous section, the vectors $v_k$ are:
$$v_1=(1,0), v_2=(1,0),v_3=(1,1), v_4=(0,1), v_5=(0,1),$$
and the drifts are:
$$\delta^1=(\lambda_1-1,\lambda_2+ \nu),\delta^2=(\lambda_1-a_1,\lambda_2+ \nu-a_2),\delta^3=(\lambda_1+\nu-a_1,\lambda_2-a_2),\delta^4=(\lambda_1+\nu,\lambda_2-1).$$

Using the results of the previous section, the interior of the stability region can be described as follows.

Assume first that $\lambda_1 -a_1 >  \nu +  \lambda_2 -a_2$.
Then one of the two following conditions should hold:
\begin{itemize}
\item
$\lambda_2 + \nu < a_2 \mbox{ and } \lambda_1 < 1+ {a_1-1 \over a_2}(\lambda_2 + \nu),$

\item or
$\lambda_2 + \nu > a_2 + \lambda_1-a_1  \mbox{ and }  \lambda_2 - a_2 < \nu.$

\end{itemize}
Symmetric conditions with the indices $1$ and $2$ reversed should hold if $\lambda_1 -a_1 + \nu >  \lambda_2 -a_2$.
\end{example}

\subsection{Bounds on the stability region of $3$ coupled processors}

Consider a process of dimension $3$ where the death rates
of each dimension depend on whether the other coordinates are strictly positive or zero,
so that for all $i\neq j\neq k$, $x_i>0$:
$$
 \phi_i(x) =
  \begin{cases}
  a_i   , \quad & x_j=0, \ x_k=0,\\
  a_{ij}, \quad & x_j>0, \ x_k=0,\\
  1,      \quad & x_j>0, \ x_k>0,
 \end{cases}
$$
which leads to the following drifts:
\begin{eqnarray*}
 \delta(x) =
 \begin{cases}
  \delta^i : ~ \delta_i^i = \lambda_i-a_i, ~  \delta_j^i=\lambda_j, ~ \delta_k^i=\lambda_k  ~, \quad \text{for} \quad x_j = 0, \quad x_k=0,\\
  \delta^{ij} : ~ \delta_i^{ij} = \lambda_i-a_{ij}, ~  \delta_j^{ij}=\lambda_j-a_{ji}, ~ \delta_k^{ij}=\lambda_k, \quad \text{for} \quad  x_j > 0, \quad x_k=0,\\
  \delta ~ = (\lambda_i-1)_{i=1\ldots 3}      \quad \text{for} \quad x_j > 0, \quad x_k>0.
 \end{cases}
\end{eqnarray*}

Let us assume $a_i \ge a_{ij} \ge 1$, so that $\phi=(\phi_1,\phi_2,\phi_3)$ is partially
decreasing.

\

Theorem 4.4.4 in \cite{fayolle1995} and
Theorem 3 in \cite{borst2008} show that the stability
region is a union of six regions corresponding to the six
possible permutations of the coordinates.
The first of these regions corresponding to the identity permutation is
the set of $(\lambda_1,\lambda_2,\lambda_3)$ such that
\begin{align}
\delta_1 &< 0,\label{c31} \\
\delta^{23}_2 &< \lambda_1(1-a_{23}) \label{c32}, \\
\delta^3_3 \pi_{00} + \delta^{13}_3 \pi_{10} + \delta^{23}_3 \pi_{01} + \delta_3 \pi_{11} &< 0 \label{c33},
\end{align}
where
\begin{align*}
\pi_{00} = P(Y_1=0,Y_2=0),~ \pi_{01} = P(Y_1=0,Y_2>0),\\
\pi_{10} = P(Y_1>0,Y_2=0), ~\pi_{11} = P(Y_1>0,Y_2>0),
\end{align*}
and $Y=(Y_1,Y_2)$ is a random vector distributed according to the stationary number of
a process in which coordinate $3$ would be always strictly positive.

\

On the other hand the sufficient conditions obtained in Section \ref{sec:discont}
can be written as the complement of the following set:

\begin{eqnarray}
\delta &> 0, \\ \label{c31b}
\mbox{or }  \exists i,j, \mbox{ and }\alpha_1,\alpha_2 \ge 0 \mbox{ such that }
\langle \alpha_1 \delta^{ij}+ \alpha_2 \delta, e_i +e_j \rangle  &>0, \\ \label{c32b}
\mbox{or } \exists i \mbox{ and } (\alpha_l)_{l=1\ldots 4} \ge 0 \mbox{ such that }
\delta^i_i \alpha_1 + \delta^{ij}_i \alpha_2 + \delta^{ik}_i \alpha_3 +\delta_i \alpha_4  &< 0 \label{c33b}.
\end{eqnarray}

It is not difficult to verify that conditions (\ref{c31}--\ref{c33}) are implied by the sufficient
conditions (\ref{c31b}--\ref{c33b}) obtained in theorem \ref{thm:discont} while the converse is not true.

\section{Conclusions} \label{sec:concl}
We derived various computable criteria of stability and instability for continuous drifts in any dimension
and for piece-wise constant drifts in dimension 2, together with generic sufficient conditions for discontinuous drifts in any dimension. An important direction of future research is to systematically characterize the second vector field and the stability conditions in dimension 3 and more.

\section{Acknowledgement}
The authors are very grateful to Mark Peletier for our fruitful discussions and his enlightening suggestions.

\end{document}